\documentclass[%
showkeys,
superscriptaddress,
nofootinbib,
preprint,
 amsmath,amssymb,
 aps,
pre,
floatfix,
]{revtex4-1}

\usepackage{graphicx}
\usepackage{dcolumn}
\usepackage{bm}
\usepackage{color}
\usepackage{comment}
\usepackage{multirow}

\let\originaleqref=\eqref
\renewcommand{\eqref}{Eq.~\originaleqref}

\begin{document}

\title{Beyond Newton: a new root-finding fixed-point iteration for nonlinear equations}

\author{Ankush Aggarwal}
\email{To whom correspondence should be addressed. E-mail: a.aggarwal@swansea.ac.uk}
\author{Sanjay Pant} 

\affiliation{Zienkiewicz Centre for Computational Engineering,\\ College of Engineering, Swansea University, Swansea SA1 8EN, UK}

\keywords{Newton-Raphson method $|$ Newton's iteration $|$ Nonlinear equations $|$ Iterative solution $|$ Gradient-based methods} 


\begin{abstract}
Finding roots of equations is at the heart of most computational science. A well-known and widely used iterative algorithm is the Newton's method. However, its convergence depends heavily on the initial guess, with poor choices often leading to slow convergence or even divergence. In this paper, we present a new class of methods that improve upon the classical Newton's method. The key idea behind the new approach is to develop a relatively simple multiplicative transformation of the original equations, which leads to a significant reduction in nonlinearities, thereby alleviating the limitations of the Newton's method. Based on this idea, we propose two novel classes of methods and present their application to several mathematical functions (real, complex, and vector). Across all examples, our numerical experiments suggest that the new methods converge for a significantly wider range of initial guesses with minimal increase in computational cost. Given the ubiquity of Newton's method, an improvement in its applicability and convergence is a significant step forward, and will reduce computation times several-folds across many disciplines. Additionally, this multiplicative transformation may improve other techniques where a linear approximation is used.
\end{abstract}


\maketitle

\clearpage

\section{Introduction}

Newton-Raphson, or often referred to as Newton's, fixed-point iteration method has been the gold-standard for numerically solving equations for several centuries. In order to set the symbols and nomenclature, we start by defining a generic problem.
\begin{equation}
\text{Find }x=x^* \text{ such that } r(x^*)=0
\label{original-eq}
\end{equation}
for a given function $r: \mathbb{K} \rightarrow \mathbb{K}$ $(\mathbb{K}=\mathbb{C}$ or $\mathbb{R})$. In general, fixed-point iteration methods start with a guess to the solution $x_n=x_0$ and iteratively update it using
\begin{equation}
x_{n+1} = \phi(x_n),\label{fixed-point-eq}
\end{equation}
where $\phi(x)$ depends on $r(x)$ and the chosen numerical scheme. It is required that $x_n\rightarrow x^*$ as $n\rightarrow\infty$ for the numerical scheme to converge, and, in practice, the fastest possible convergence is desired. If we expand $r(x)$ in finite Taylor's series up to second order about point $x_n$ and evaluate it at $x=x^*$ to solve $r(x^*)=0$, we get
\begin{equation}
r(x^*)=r(x_n) + \left.\frac{dr}{dx} \right|_{x_n} (x^*-x_n) + \left.\frac{1}{2}\frac{d^2 r}{d x^2} \right|_{\zeta} (x^*-x_n)^2=0,
\label{taylor-expansion}
\end{equation}
for some $\zeta \in [x_n,x^*]$. The last term is called the remainder and can be written in other forms \cite{apostol1967calculus}. Here, $\zeta$ is unknown. Therefore, in the standard Newton's method, we neglect the second order term and approximate $x^*$ with an updated guess to the solution: $x_{n+1}$. Thus, we obtain the Newton's fixed-point iteration, also called Newton-Raphson method (using notation $r'=dr/dx$):
\begin{equation}
x_{n+1}  =  x_n + \Delta x^{\textrm{N}}   = x_n - \frac{r(x_n)}{r'(x_n)}.
\label{original-newton}
\end{equation}
Using \eqref{taylor-expansion} and \eqref{original-newton}, one can rearrange to get the evolution of error $e_n=|x^*-x_n|$:
\begin{equation}
e_{n+1} = \left|\frac{r''(\zeta)}{2r'(x_n)} \right| e_n^2.
\label{error-evolution}
\end{equation}

The above step shows (at least) quadratic convergence of the Newton's method when $r(x)$ has only simple roots. This quadratic convergence has led to a wide adoption of Newton's method in solving problems in all scientific disciplines. Furthermore, this attractive property has also led to a large amount of work towards developing root-finding algorithms that either imitate or approximate the Newton's method.

\begin{figure}[!h]
\centering
\includegraphics[width=0.55\columnwidth]{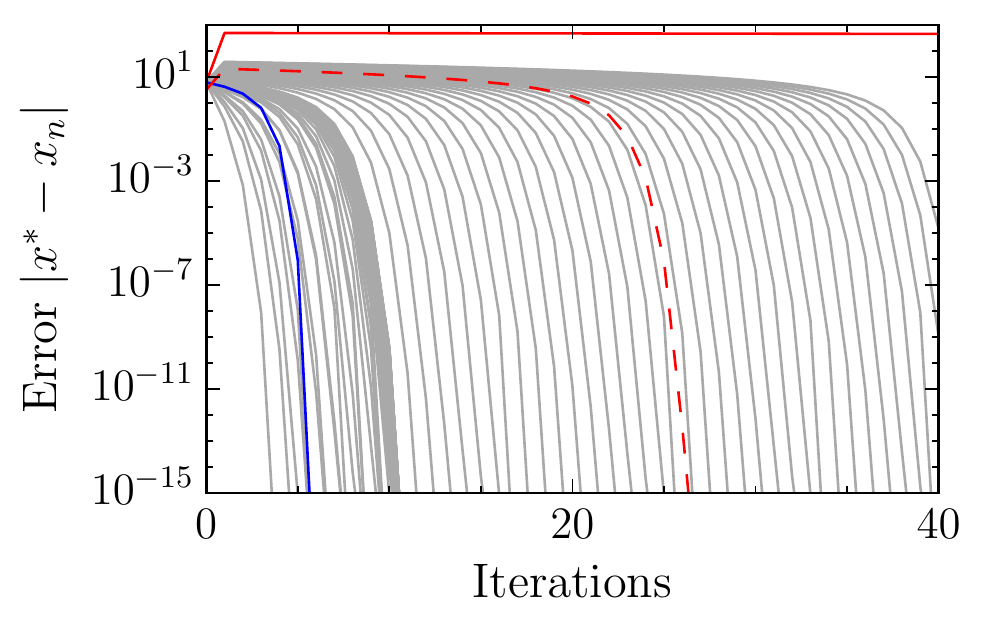}
\caption{Convergence of the standard Newton's method (red lines), the Extended Newton method (grey lines), and the Corrected Newton method (blue line) for solving $e^x-500=0$ (solid curves with initial guess $x_0=0$ and dashed curve with initial guess $x_0=3$)}
\label{ex1}
\end{figure}

Despite the aforementioned attractive convergence of the Newton's method, it can be seen from \eqref{error-evolution} that to achieve quadratic convergence, the current guess $x_n$ must be ``sufficiently close'' to the solution. When the initial guess is not close to the solution, the convergence can be slower than quadratic, or, worse, the iterations may not converge, oscillate, or diverge. As an example, if we consider $r(x)=e^x-500$ and solve it using \eqref{original-newton} with $x_0=0$, iterations diverge until we get numerical overflow (Fig.~\ref{ex1}, solid red line). Instead, if we start with $x_0=3$, the error decreases sub-quadratically before achieving quadratic convergence closer to the solution (Fig.~\ref{ex1}, dashed red line).

\section{Methods and Results}
The basin of attraction for a root $x^*$ is defined as the set of initial guesses from which a given fixed-point iteration converges to $x^*$ (this can be seen as the domain of convergence for initial guesses). Naturally, it is desirable to have a large basin of attraction for the root so that convergence is achieved for a wide range of initial guesses. Even though, for a general $r(x)$, determination of the size of the basin of attraction is challenging, it is clear from \eqref{error-evolution} that size of the basin depends on a measure of nonlinearity\footnote{This nonlinearity measure is similar to the $\gamma$ number in Smale's $\alpha$-theory\cite{smale2000algorithms} } $N(\zeta,x_n):=\left|\dfrac{r''(\zeta)}{2r'(x_n)} \right|$.

We pose a question: can we pre-multiply the original equation by a function so that we obtain a larger (or, if possible, infinite) basin of attraction? That is, instead of solving \eqref{original-eq}, can we solve
\begin{equation}
\text{Find }x=x^* \text{ such that } P(x^*)r(x^*)=0,
\label{new-eq}
\end{equation}
\emph{still using Newton's method}? The idea is to choose a $P(x)$ that decreases the nonlinearity $N$ and, hence, gives a larger basin of attraction, while retaining at least quadratic convergence close to the root. For \eqref{new-eq}, the measure of nonlinearity is 
\begin{equation}
N(\zeta,x_n) = \left| \frac{\left[ P(x)r(x) \right]''|_{x=\zeta}} {\left[ P(x)r(x) \right]'|_{x=x_n}} \right|.
\label{nonlinearity-measure}
\end{equation}
We note that if $P(x)$ is not a function of $x$, we do not change the nonlinearity of the problem, and this case would fall under the purview of the highly developed field of linear preconditioning. In an attempt to minimize $N$, \eqref{nonlinearity-measure}, we equate $N=0$, which gives 
\begin{equation}
P(x) = \frac{c_1x-c_2}{r(x)},\label{eqn_px_orig}
\end{equation}
for arbitrary integration constants $c_1$ and $c_2$. However, this has two problems. Multiplying the above form of $P(x)$ by $r(x)$: 
\begin{enumerate}
\item eliminates the original root $x=x^*$ and
\item introduces a new root $x=c_2/c_1$.
\end{enumerate}
In order to to solve these problems, we add another constant $\kappa\ne0$ in the denominator of $P(x)$ in \eqref{eqn_px_orig}
\begin{equation}
P(x) = \frac{c_1x-c_2}{r(x)+\kappa}.
\end{equation}
This avoids the elimination of the $x^*$ root. In addition, we eliminate the undesirable new root $x=c_2/c_1$ by choosing $\kappa$ such that $r(c_2/c_1)+\kappa=0$. That is, $\kappa=-r(c_2/c_1)$. Furthermore, $c_1=0$ makes $c_2/c_1$ indeterminate and any value $c_1\ne0$ only scales the equation with a constant without affecting its nonlinearity. Therefore, without loss of generality, we set $c_1=1$ and $c_2=c$. That is,
\begin{equation}
P(x) = \frac{x-c}{r(x)-r(c)}.
\end{equation}
Thus instead of applying Newton's method to $r(x)=0$, if we apply Newton's method to our new equation
\begin{equation}
\frac{(x-c)r(x)}{r(x)-r(c)}=0,
\label{new-equation}
\end{equation}
we get a new fixed-point iteration
\begin{equation}
x_{n+1} = x_n + \Delta x^{\textrm{EN}}  = x_n  - \frac{(x_n-c)r(x_n)}{r(x_n)-(x_n-c)r'(x_n)\frac{r(c)}{r(x_n)-r(c)}}.
\label{new-newton}
\end{equation}
The superscript $^{\textrm{EN}}$ denotes ``Extended Newton'' method.  

\subsection*{Choice of $c$}

The new fixed-point iteration \eqref{new-newton} contains an arbitrary constant $c$, which remains to be determined. We note that if $c=x^*$ such that $r(c)=0$, we have $x_{n+1} = x_n + \Delta x^{\textrm{EN}}  = x_n +( c - x_n ) = c = x^*$. That is, we will find the solution in  a single iteration irrespective of the starting point $x_0$ and function $r$. This is not surprising; if one could chose $c=x^*$, the solution is already known, making its basin of attraction infinite.

However, it is not clear how the distance of $c$ from the solution will affect convergence. Even though one could study the effect of $c$ on the nonlinearity by assessing $N(\zeta,x_n)$ from \eqref{nonlinearity-measure}, it does not give an immediate insight. For the $r(x)=e^x-500$ example, the new equation to be solved is $(x-c)(e^x-500)/(e^x-e^c)=0$. Surprisingly, we find that our Extended Newton method, \eqref{new-newton}, converges starting from $x_0=0$ (Fig.~\ref{ex1}, gray lines), irrespective of the choice of $c$ in a wide range $c\in(-50,50)$. 
Even if this behaviour of $c$ is specific to the $r(x)$ chosen, given the simplicity of the formulation, this is a remarkable result showing the potential of the EN method in comparison to the classical Newton's method. 


\subsection*{Limiting case of $c\rightarrow x_n$}
Next, we note that we must choose $c\ne x_n$. Choosing $c=x_n$ makes $x_n$ a stationary point of the fixed-point iteration \eqref{new-newton}, and the iteration gets ``stuck'' at this undesired point. This is because, even though $\lim_{x\rightarrow c} P(x)r(x) \ne 0$, $P(c)r(c)$ is indeterminate. Practically, this does not pose a problem as long as $c\ne x_0$. 

Nevertheless, we look at the limiting case when $c\rightarrow x_n$. Rewriting \eqref{new-newton} and taking the limit we get
\begin{equation}
\lim_{c\rightarrow x_n} \Delta x^{\textrm{EN}}  = -\frac{r(x_n)/r'(x_n)}{1-r(x_n)r''(x_n)/2r'(x_n)^2}.
\label{limit-c-deltax}
\end{equation}
This gives us a new fixed-point iteration
\begin{equation}
x_{n+1} = x_n + \Delta x^{\textrm{CN}}  = x_n -\frac{r(x_n)/r'(x_n)}{1-r(x_n)r''(x_n)/2r'(x_n)^2},
\label{limit-c-deltax}
\end{equation}
where the superscript $^{\textrm{CN}}$ denotes ``Corrected Newton'' method (the choice of this name will become clear shortly). The advantage of this fixed-point iteration is that there is no arbitrary constant involved, however we have to pay with the price of calculating the second derivative of $r(x)$
\footnote{The same limit on our modified equation \eqref{new-equation} gives
\begin{equation}
\lim_{c\rightarrow x} \frac{(x-c)r(x)}{r(x)-r(c)}=\frac{r(x)}{r'(x)}=0.
\end{equation}
If we apply Newton's method to the above equation, we get
\begin{equation}
\Delta x  = -\frac{r(x_n)/r'(x_n)}{1-r(x_n)r''(x_n)/r'(x_n)^2}.
\label{rbyrprime-deltax}
\end{equation}
There is a factor of 2 between $\Delta x^{\textrm{CN}}$ and \eqref{rbyrprime-deltax}. Numerically, we verified that \eqref{limit-c-deltax} performs better than \eqref{rbyrprime-deltax}.}. 

To analyze \eqref{limit-c-deltax} further, we re-write the Taylor's expansion \eqref{taylor-expansion} with $\zeta=x_n$ (thereby neglecting the remainder term of third order) and $\Delta x=x_{n+1}-x_n$ as
\begin{equation}
r(x_n) + \left( r'(x_n)+ \frac{r''(x_n)}{2}\Delta x\right) \Delta x =0.
\end{equation}
If we substitute $\Delta x$ within the parentheses as the one from standard Newton \eqref{original-newton} and rearrange, we get \eqref{limit-c-deltax}. That is, \eqref{limit-c-deltax} can be thought of as a correction to the Newton's equation using second derivative $r''(x)$. Thus, we can re-write \eqref{limit-c-deltax} as 
\begin{equation}
r(x_n) + \left( r'(x_n)+ \frac{r''(x_n)}{2}\Delta x^{\textrm{N}} \right) \Delta x^{\textrm{CN}}  =0,
\label{correction-equation}
\end{equation}
where $\Delta x^{\textrm{N}} $ is obtained from standard Newton's method. While the above form is obtained with the aim of expanding the basin of attraction, we note that the \eqref{correction-equation} is identical to the Halley's method with cubic order convergence close to the root \cite{traub1982iterative}. This observation will be useful in developing the extension of the method to vector functions.

\begin{figure}[!h]
\centering
\includegraphics[width=0.7\columnwidth]{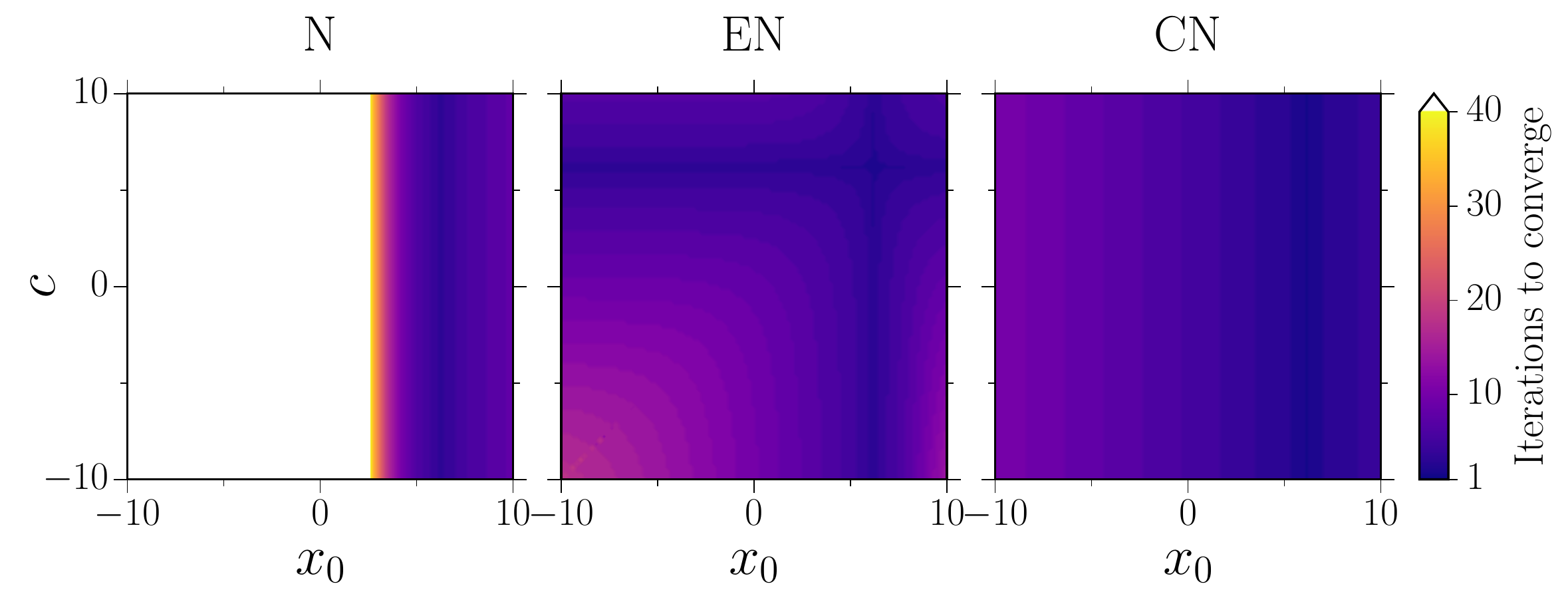}
\includegraphics[width=0.7\columnwidth]{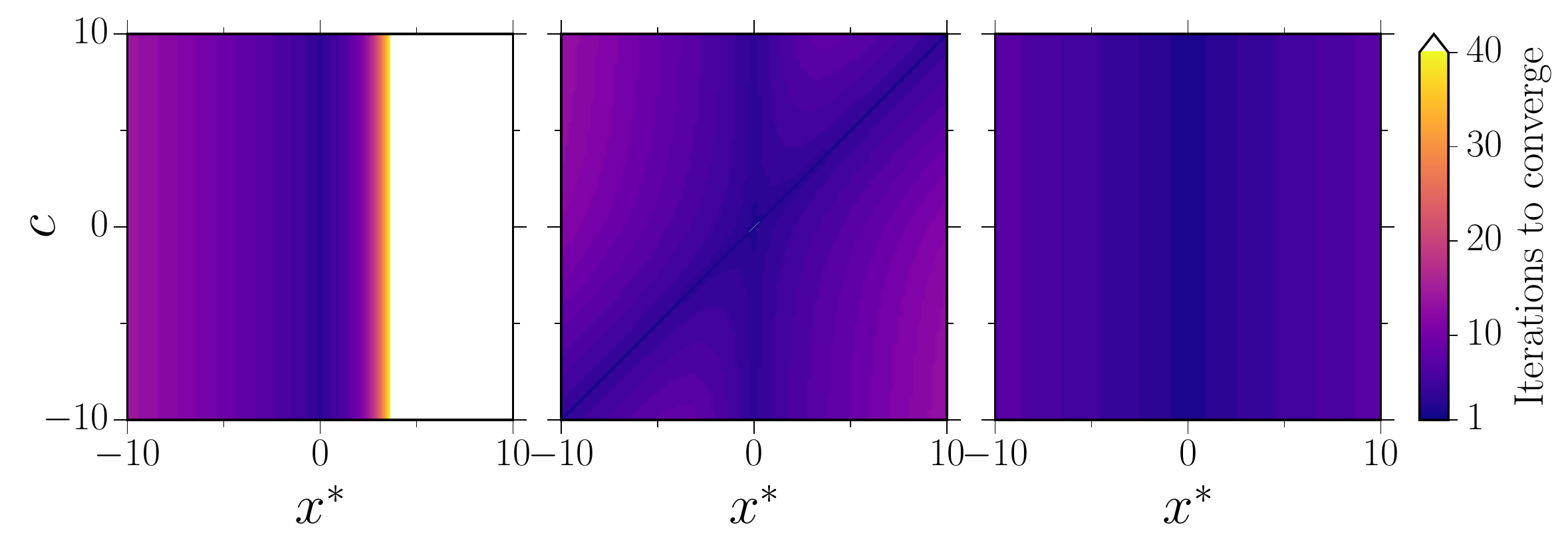}
\caption{(From left) Iterations taken to converge for the standard Newton's method, the Extended Newton method, and the Corrected Newton method for solving $r(x)=e^x-H=0$: with $H=500$ and varying initial guess $x_0$ (top), and with $x_0=0$ and varying $H=e^{x^*}$ (bottom)}
\label{ex2}
\end{figure}

We solve the problem $e^x-H=0$ for varying values of $c$, $x_0$, and $H$ and find a significantly improved convergence using the proposed fixed-point iterations (Fig.~\ref{ex2}). In order to verify the generality of the new schemes, we solve several different nonlinear equations and provide the results in the Supporting information (SI).

\subsection*{Complex Plane}
\begin{figure}[!h]
\centering
\includegraphics[width=0.7\columnwidth]{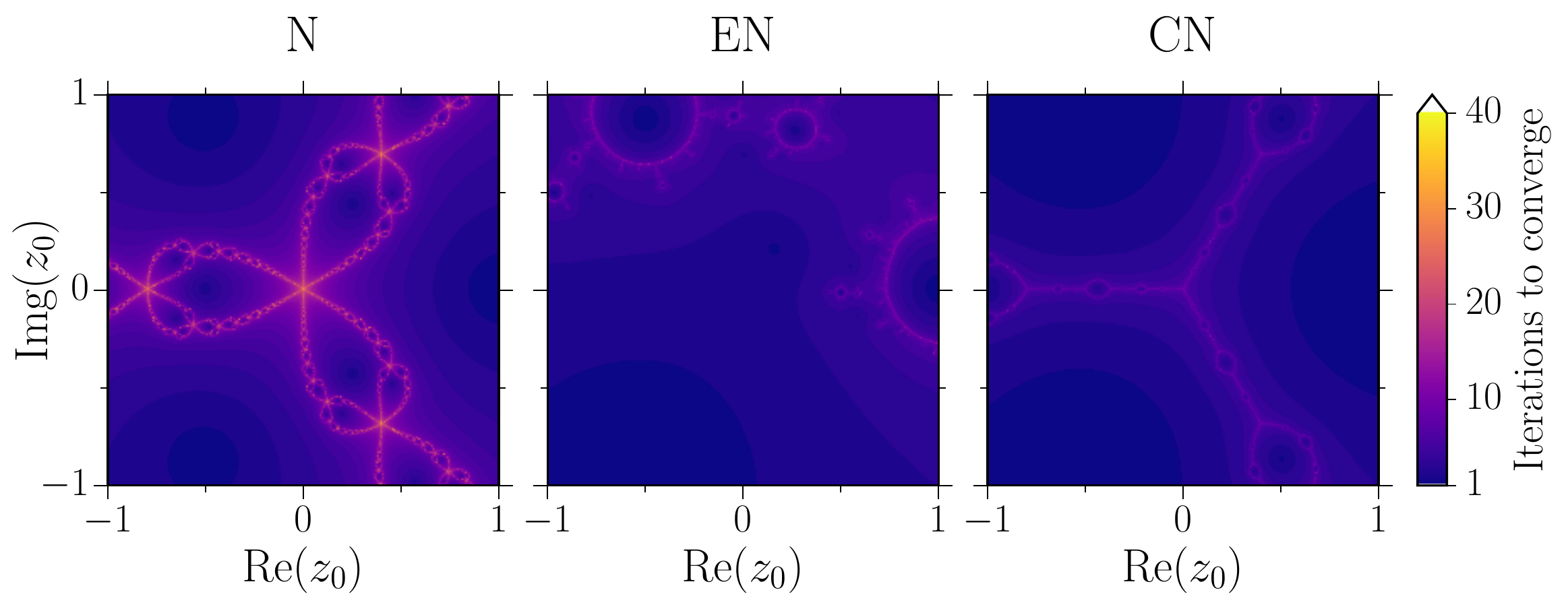}
\includegraphics[width=0.7\columnwidth]{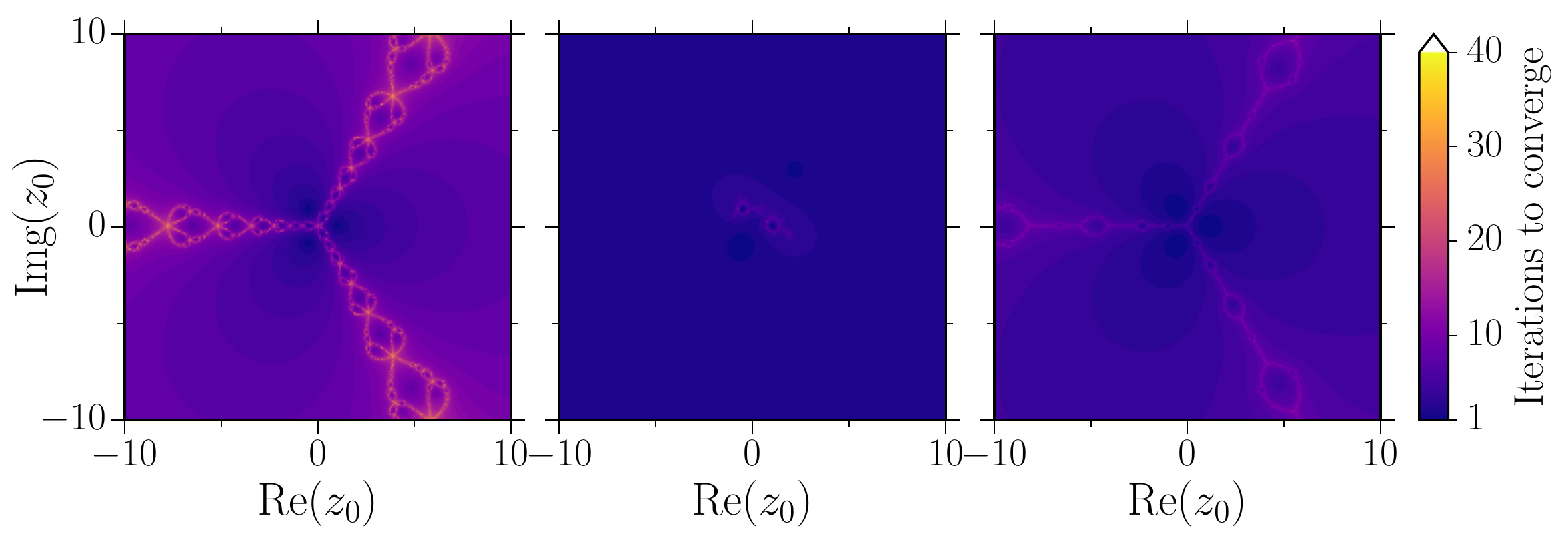}
\includegraphics[width=0.7\columnwidth]{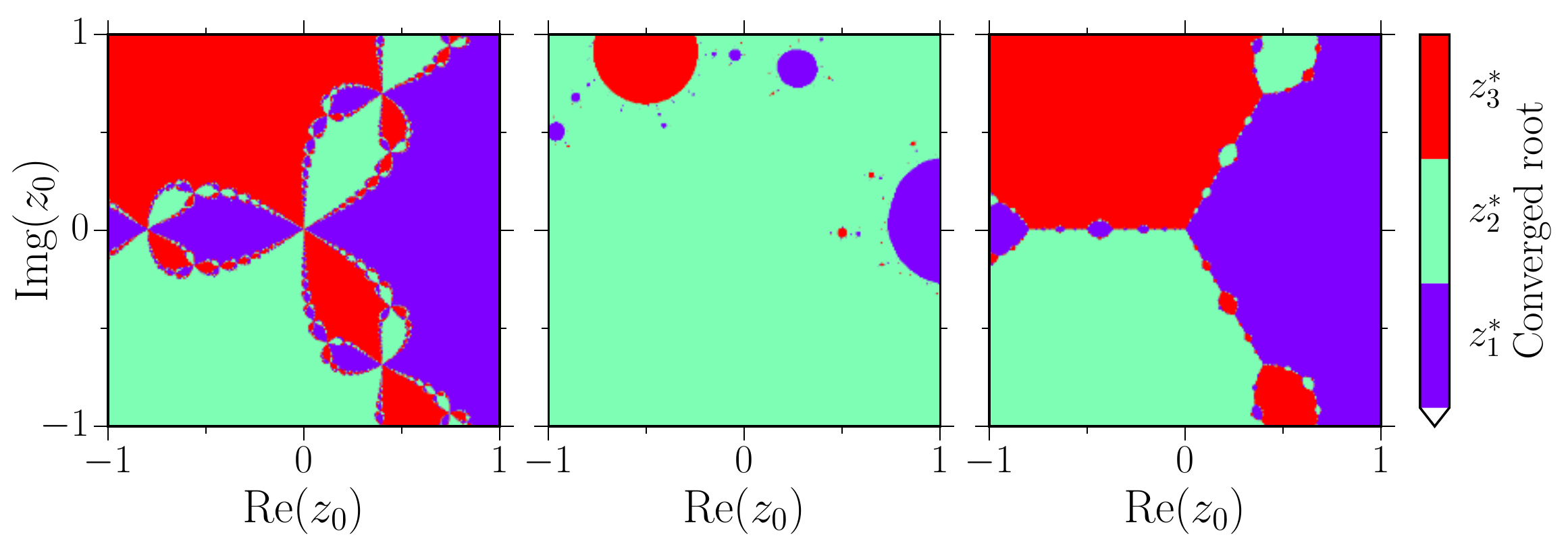}
\includegraphics[width=0.7\columnwidth]{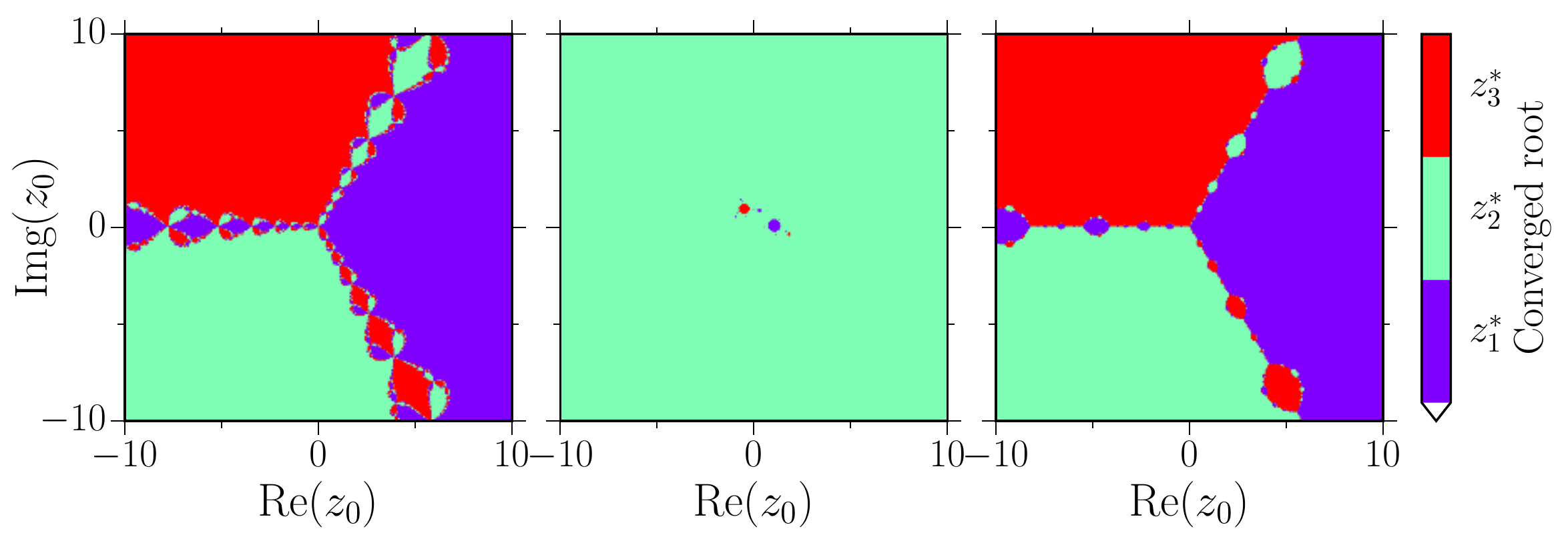}
\caption{Newton fractal for $r(z)=z^3-1$ in complex plane, (top two rows) colored by number of iterations and (bottom two rows) colored by the final root:  the standard, Extended, and Corrected Newton methods (from left to right). $c=-0.65-0.65i$ is used for EN, and ${z}^*_1=1+0i$, ${z}^*_2=-0.5-\sqrt{0.75}i$ and ${z}^*_3=-0.5+\sqrt{0.75}i$.}
\label{newton-fractal}
\end{figure}

The extension to functions on the complex plane is straightforward, and the same equations can be applied. Applying Newton's method to $r(z)=z^3-1$ gives us a Newton fractal, and the Extended Newton and Corrected Newton methods modify that fractal (Fig.~\ref{newton-fractal}). Surprisingly, the new schemes even reduce the dimension of the fractal (Fig.~S8) and the basins of attraction appear more connected (generally speaking) than before. Choosing $c$ in Extended Newton breaks the three-fold symmetry of this system, and the solution close to $c$ is heavily favored.

\subsection*{Extension to multiple unknowns}
Instead of single unknown, we have $m$ unknowns $x_i$ and a set of $m$ equations $r_i(x_1,\dots,x_m)=0$ for $i=1,\dots, m$. Our problem in $m$-dimensions is
\begin{equation}
\text{Find }x_i=x_i^* \text{ such that } r_i(\mathbf{x}^*)=0 \;\forall i=1,\dots,m,
\label{original-multiple-eq}
\end{equation}
where all $m$ unknowns are collectively written as $\mathbf{x}$. Applying the standard Newton's method gives us a system of $m$ linear equations to be solved to obtain the step $\Delta x_j^{\textrm{N}}$
\begin{equation}
\sum_j r_{i,j}(\mathbf{x}) \Delta x_j^{\textrm{N}}  = -r_i(\mathbf{x}),
\label{original-newton-multiple}
\end{equation}
where $r_{i,j}=\partial r_i/ \partial x_j$ and we have dropped the subscript $n$ to denote the current iteration value of $\mathbf{x}$ as that is implied.

There could be several different way of extending \eqref{new-equation} to this case of multiple unknowns. If it is possible to recognize a subset of the equations that are most nonlinear, one may use Extended Newton method for that subset of equations and standard Newton's method for the rest. Here, we consider only the case where all equations are transformed. Thus, one straightforward extension is the following modified set of equations:
\begin{equation}
\frac{(x_i-c_i)r_i(\mathbf{x})}{r_i(\mathbf{x})-r_i(\mathbf{x}_{ci})}=0
\end{equation}
for chosen values of $c_i$ ($i=1,\dots,m$), and $\mathbf{x}_{ci}$ is the vector $\mathbf{x}$ with $i$-th element replaced by $c_i$. Accordingly, applying the Newton's method we get our new system of linear equations 
\begin{align}
\sum_j\Bigg\{\delta_{ij}\frac{r_i(\mathbf{x})}{r_i(\mathbf{x})-r_i(\mathbf{x}_{ci})} + \frac{(x_i-c_i)}{\left[r_i(\mathbf{x})-r_i(\mathbf{x}_{ci})\right]^2}  \nonumber \\
 \big[ r_{i,j}(\mathbf{x}_{ci}) r_i(\mathbf{x}) (1-\delta_{ij}) - r_{i,j}(\mathbf{x}) r_i (\mathbf{x}_{ci})  \big]  \Bigg\}\Delta x_j^{\textrm{EN}} \nonumber \\
  = -\frac{(x_i - c_i) r_i(\mathbf{x})}{r_i(\mathbf{x})-r_i(\mathbf{x}_{ci})}, \label{extended-newton-multiple}
  \end{align}
where $\delta_{ij}=1$ when $i=j$ and 0 otherwise. For real or complex functions with single unknown, we have to calculate $r(c)$ only once and there is no need for computing $r'(c)$. However, in the above extension to multiple unknowns, $r_i(\mathbf{x}_{ci})$ needs to be evaluated at every iteration along with $r_{i,j}(\mathbf{x}_{ci})$ as well. Thus, the Extended Newton method for the case with multiple unknowns leads to a cumbersome set of calculations. 

On the other hand, extending the Corrected Newton method, \eqref{correction-equation}, to multiple unknowns is simpler:
\begin{equation}
\sum_j \left( r_{i,j} + \frac{1}{2} \sum_k r_{i,jk} \Delta x_k^{\textrm{N}} \right) \Delta x_j^{\textrm{CN}} = -r_i,
\label{correction-equation-multiple}
\end{equation}
where $\Delta x_j^{\textrm{N}}$ is calculated from \eqref{original-newton-multiple} and $r_{i,jk}=\dfrac{\partial^2 r_i}{ \partial x_j \partial x_k}$. Thus, \eqref{original-newton-multiple} and \eqref{correction-equation-multiple} can be solved successively.

The above method requires solution to two systems of $m$ linear equations. To alleviate this and obtain only one system of $m$ equations, we multiply both sides of \eqref{correction-equation-multiple} by $r_{i,i}$
\begin{equation}
\sum_j \left( r_{i,i}r_{i,j} + \frac{1}{2} \sum_k r_{i,i}r_{i,jk} \Delta x_k^{\textrm{N}} \right) \Delta x_j^{\textrm{CN}} = -r_ir_{i,i}.
\end{equation}
Using an approximation that $r_{i,i}r_{i,jk} \approx r_{i,k} r_{i,ji}$ and utilizing \eqref{original-newton-multiple}, we get a simplified expression
\begin{equation}
\sum_j \left( r_{i,i}r_{i,j} - \frac{1}{2} r_{i,ji} r_i \right) \Delta x_j^{\textrm{QCN}} = -r_ir_{i,i}.
\label{correction-equation-multiple2}
\end{equation}
This equation gives us the step $\Delta x_j^{\textrm{QCN}}$ by solving only one system of $m$ equations, and we call it a quasi-Corrected Newton method (QCN). In numerical methods, computing the Jacobian $r_{i,j}$ is common, but calculating another gradient $r_{i,jk}$ for \eqref{correction-equation-multiple} may not be desirable. However, in \eqref{correction-equation-multiple2}, we note that if we have an expression for $r_{i,j}$, $r_{i,ji}$ can be computed using finite difference in only $\mathcal{O}(m)$ operations, making this approach highly attractive.

\begin{figure}[!h]
\centering
\includegraphics[width=0.55\columnwidth]{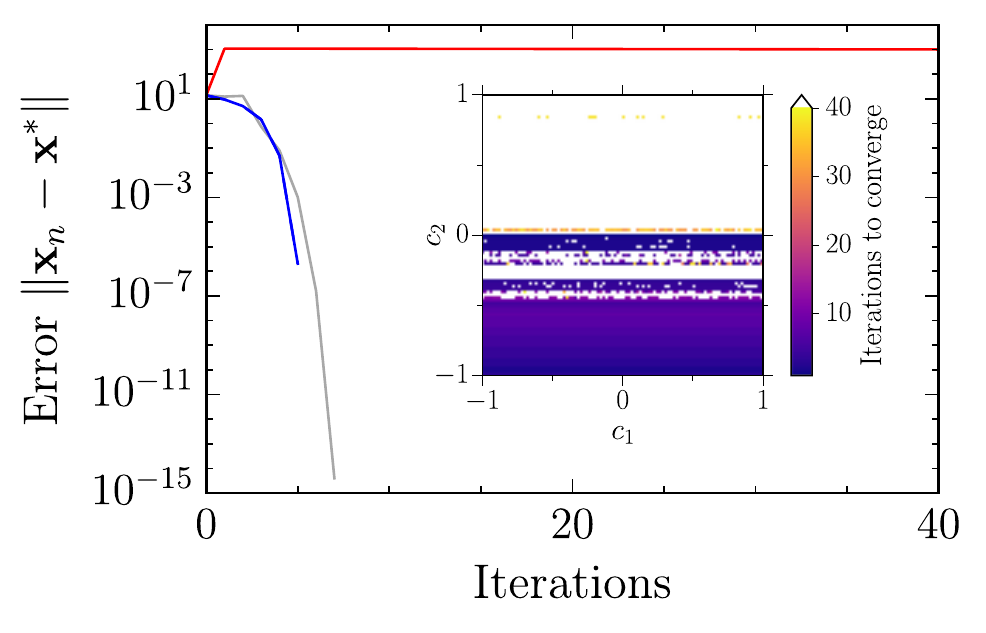}
\caption{Convergence for solving \eqref{exp-two-var} with initial guess $x_1=x_2=0$ using the standard Newton's method (red line), the Extended Newton method with $c_1=c_2=-0.5$ (grey line), and the Corrected Newton method (blue line). The inset shows convergence of EN with the same initial condition and varying $c_1$ and $c_2$.}
\label{multi-ex}
\end{figure}

As an extension of the exponential example from single variable system, we solve the following two equations:
\begin{align}
e^{x_1}-e^{x_2-x_1} &= 0 \nonumber \\
e^{x_2-x_1}-1 &= 500.
\label{exp-two-var}
\end{align}
This represents two springs in series with one end fixed and the other end under force $H$, where each spring's force is nonlinearly related to its change in length as $e^{\Delta l}-1$. Similar to the single equation, the standard Newton's method fails to converge for initial guess $x_1=x_2=0$, whereas the Extended Newton (for $c_2<-0.5$) and Corrected Newton methods converge in less than 10 iterations (Fig.~\ref{multi-ex}).

\begin{figure}[!h]
\centering
\includegraphics[width=0.7\columnwidth]{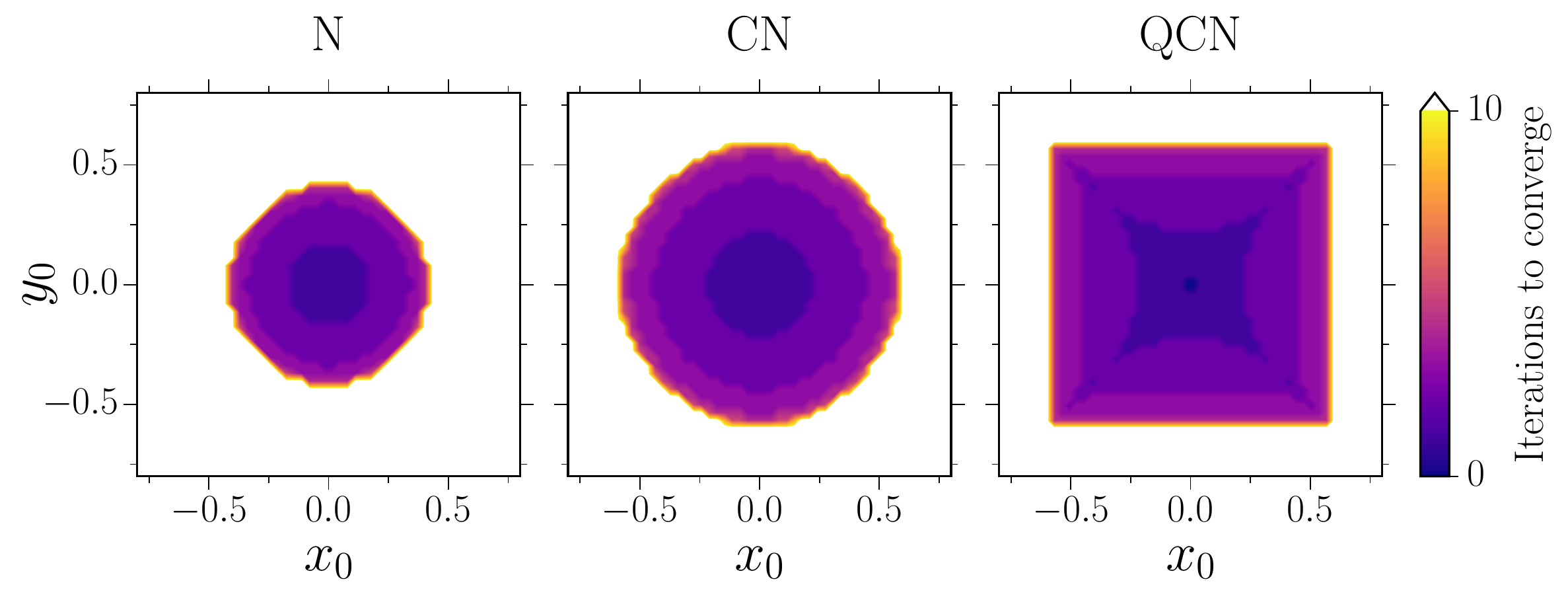}
\includegraphics[width=0.7\columnwidth]{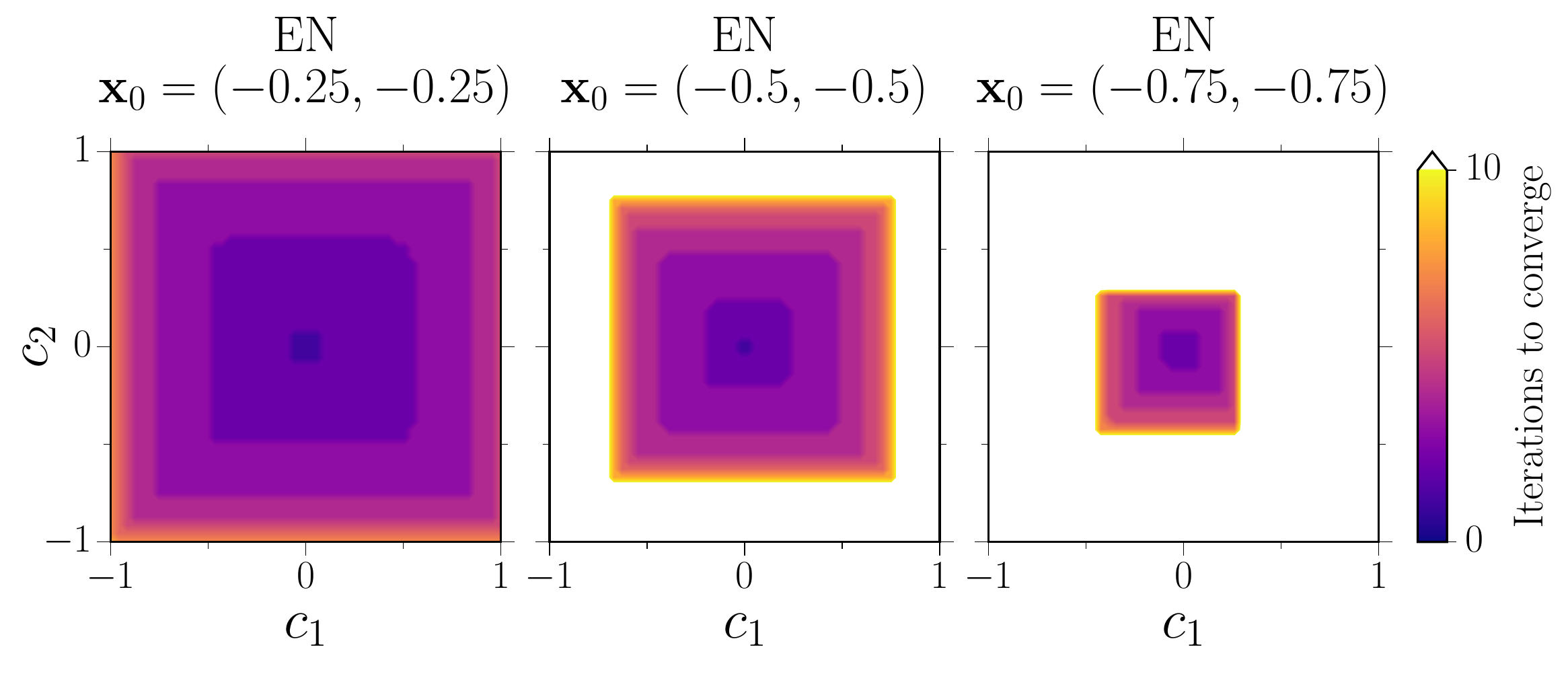}
\caption{Convergence in finding the minimum of two-variable Easom's function \eqref{easom-func}: (top) for the standard and Corrected Newton method with varying initial conditions, and (bottom) for Extended Newton method with varying $\mathbf{c}$ and three initial conditions.}
\label{easom-result}
\end{figure}

As another example, we consider the minimization of Easom's function with two variables:
\begin{equation}
f(x,y)=-\cos(x)\cos(y)\exp(-x^2-y^2).\label{easom-func}
\end{equation}
For this function with a single minimum at $(0,0)$, the basin of attraction is relatively small using the standard Newton's method. Our proposed methods increase that basin significantly (Fig.~\ref{easom-result}). Surprisingly, the improvement using quasi-Corrected Newton method is even greater than the Corrected Newton method. Furthermore, the Extended Newton method can provide an even larger basin of attraction if $c$ is chosen close enough to the solution (Fig.~\ref{easom-result} bottom row). A summary of proposed methods is provided in Table~1.
%


\section*{Discussion}
Newton's method is one of the most widely used methods for finding roots in science and engineering. This is particularly true when the computation of the Jacobian is inexpensive. In the opposite case, significant amount of work has focused on developing methods that approximate the Newton's method, largely through approximations of $r_{i,j}$ leading to the so-called quasi-Newton methods, so that convergence properties similar to the Newton's method can be achieved with reduced computational expense \cite{morton1995iterative}. On the other hand, different methods have been proposed that achieve higher than quadratic convergence rates using higher order derivatives. For example, Halley's method, Cauchy's method, Chebyshev method \cite{traub1982iterative}, class of Householder method \cite{householder1970numerical}, and other variants \cite{osada1994optimal, amat2003geometric, kou2008some, cordero2013dynamics}. However, the convergence of all of these methods is examined when the starting guess is close to the solution. 

The objective of this article is to go beyond the classical approaches. The central idea is to decrease the system nonlinearity in view of improving convergence properties \emph{far from the solution}. Starting with one variable and one equation, we premultiplied the original equation and formulated a new equation with the same roots (Eq.~\ref{new-equation}). Applying Newton's method to this equation, we proposed a new fixed-point iteration, which we term as the Extended Newton (EN) method \eqref{new-equation}. 
%
%
This method has an arbitrary constant $c$, and, in fact, the choice of $c$ makes Extended Newton method a family of fixed-point iterations. Using several example functions, we found that the Extended Newton method, unlike standard Newton's method, can find the solution even when the starting point is considerably far from the solution. We found that irrespective of the choice of $c$, the performance of EN was always better than the standard Newton's method. It is also important to recognize that even if we choose $c$ close to the starting guess, the convergence is greatly improved. Therefore, when no prior information on the choice of $c$ is available, we recommend choosing $c=x_0+\epsilon$ for small $\epsilon$.

\begin{table}[h!]
\begin{tabular}{ccp{7cm}}
\hline
Method & Scalar function update & Vector function update \\
\hline
N & $\Delta x^{\textrm{N}}=- \frac{r(x_n)}{r'(x_n)}$ & $\sum_j r_{i,j}(\mathbf{x}) \Delta x_j^{\textrm{N}}  = -r_i(\mathbf{x})$ \\[10pt]
EN & $\Delta x^{\textrm{EN}}  =  - \frac{(x_n-c)r(x_n)}{r(x_n)-(x_n-c)r'(x_n)\frac{r(c)}{r(x_n)-r(c)}}$ & 
$\sum_j\Bigg\{\delta_{ij}\frac{r_i(\mathbf{x})}{r_i(\mathbf{x})-r_i(\mathbf{x}_{ci})} + \frac{(x_i-c_i)}{\left[r_i(\mathbf{x})-r_i(\mathbf{x}_{ci})\right]^2}  
 \big[ r_{i,j}(\mathbf{x}_{ci}) r_i(\mathbf{x}) (1-\delta_{ij}) - r_{i,j}(\mathbf{x}) r_i (\mathbf{x}_{ci})  \big]  \Bigg\}\Delta x_j^{\textrm{EN}} = -\frac{(x_i - c_i) r_i(\mathbf{x})}{r_i(\mathbf{x})-r_i(\mathbf{x}_{ci})}$ \\[20pt]
CN & $\Delta x^{\textrm{CN}}  =  -\frac{r(x_n)/r'(x_n)}{1-r(x_n)r''(x_n)/2r'(x_n)^2}$ & $\sum_j \left( r_{i,j} + \frac{1}{2} \sum_k r_{i,jk} \Delta x_k^{\textrm{N}} \right) \Delta x_j^{\textrm{CN}} = -r_i$ \\[20pt]
QCN & $-$ &  $\sum_j \left( r_{i,i}r_{i,j} - \frac{1}{2} r_{i,ji} r_i \right) \Delta x_j^{\textrm{QCN}} = -r_ir_{i,i}$ \\
\hline
\end{tabular}
\caption{Summary of Newton's method and its proposed extensions}
\end{table}

We note that although we used a constant value of $c$ during iteration, there are other options of choosing its value at every iteration based on the current guess $x_n$, previous guess $x_{n-1}$, and/or the function form $r(x)$. In fact, if we choose $c=x_{n-1}$, the EN can be thought of as a combination of secant method and Newton's method. As a limit case $c\rightarrow x_n$, we rediscovered Halley's method \eqref{correction-equation}, which we call the Corrected Newton (CN) method. While Halley's method was motivated by achieving cubic convergence, our approach has been motivated by increase in the size of the basin of attraction. It is, therefore, remarkable that both approaches lead to identical algorithm. Thus, one contribution of this work is to elucidate that Halley's method not only provides cubic convergence but also results in a larger basin of attraction. We also note that the Extended Newton method provides at least quadratic convergence close to the root, because it is a Newton's method applied to $P(x)r(x)$. In this sense, the EN method is no worse that Newton's method while providing a larger basin of attraction. The advantage of the Corrected Newton over Extended Newton method is that there is no arbitrary constant to choose; however, it requires the calculation of the second derivative $r''(x)$. With symbolic computation \cite{buchberger1996symbolic} and automatic differentiation \cite{rall1981automatic} becoming commonplace, this limitation is becoming less important. 

The application of these new methods to a sample of problems provided incredible improvement. The equation $e^x-H=0$ could be solved for extremely high $H$ in less than 10 iterations. Similar improvements were found for other examples (presented in the SI). We note that for a specific (known) equation form, instead of pre-multiplying, a transformation may be more suitable to decrease its nonlinearity. For example, for equation $e^x-H=0$, one may write the equation as $e^x=H$ and then take a logarithm on both sides resulting in a linear equation that will always be solved in one iteration. However, this approach is problem-specific, requires identifying the two sides of the equation, and has been proposed by us for nonlinear elasticity \cite{2018arXiv180106016M}. On the other hand, the present approach is completely general and works for all kinds of equations and nonlinearities.

The extension of the approach to functions on the complex plane was straightforward and produced some astonishing results when applied to Newton fractal. The new method not only finds the root in fewer iterations, it also decreases the sensitivity of which root is obtained with respect to perturbations in the initial guess. This was especially true for Extended Newton method, where a choice of $c$ made almost all initial guesses converge to the same solution closest to $c$. Corrected Newton method did not have such a bias and reduced the sensitivity equally for all the three roots. 

Further extension of the Extended Newton method to arbitrary number of unknowns and equations can be accomplished in multiple ways. We presented one option (Eq.~\ref{extended-newton-multiple}) and found the results improved significantly. This approach requires choosing $m$ arbitrary unknowns $c_i$ and calculation of the residual $r_i$ and Jacobian $r_{i,j}$ at $m$ additional points. This can lead to additional computation for large systems. Instead, one could choose a smaller subset of equations responsible for the nonlinearity and only transform those. This would decrease the extra computation required while still improving the convergence. Furthermore, for the exponential example, it was observed that the method worked only for $c_2<-0.5$. The underlying reason for this remains unclear, especially because the application of this method to Easom's function did not show any such behavior.

On the other hand, the extension of Corrected Newton method to system of equations was straightforward and unique (Eq.~\ref{correction-equation-multiple}). Even though this approach worked perfectly, it has two computational disadvantages: calculation and storage of $r_{i,jk}$ with $m^3$ terms and solving two linear systems -- first for the standard Newton (Eq.~\ref{original-newton-multiple}) and then for the correction (Eq.~\ref{correction-equation-multiple}). As a simplification, we proposed its approximation (Eq.~\ref{correction-equation-multiple2}), which we call quasi-Corrected Newton method. This approach only requires calculation of one diagonal of $r_{i,jk}$, i.e. $r_{i,ij}$ and solving only one linear system of equations. The name quasi-Corrected Newton method is inspired from multitude of quasi-Newton algorithms that are used to imitate Newton's method while decreasing the computational requirements. We note that even if explicit expressions for $r_{i,jk}$ may be unavailable or cumbersome to calculate, $r_{i,ij}$ can be calculated using finite difference in only $\mathcal{O}(m)$ operations, and can be easily parallelized (unlike the iterations of Newton's method).

To our knowledge, the quasi-Corrected Newton method has no hitherto reported counterpart. We believe that one of the reasons why these higher-order methods (Halley's method, Householder's, or other higher-order methods) are not widely used in science and engineering is that the faster convergence rates of these methods are usually offset by their higher computational cost at each iteration. In this respect, the quasi-Corrected Newton method provides a good balance between improved convergence properties and computational expense. In addition, the larger basin of attraction also outweighs the small extra computational cost of these proposed methods.

Among the different methods proposed, EN is associated with the least computational expense, and in many cases outperforms the Corrected Newton method, provided a good values of $c$ is available (see SI). We believe that in the future, different scientific communities will report optimal ranges of $c$, obtained through numerical analysis and/or experimentation, for specific problems of wide interest. We foresee a large impact of the proposed methods in numerical solution of partial differential equations, where the nonlinearities only allow a step-by-step increase in boundary conditions using the classical Newton's method post discretization (a common example being the nonlinear finite element analysis in structural mechanics where the applied loads can only be increased incrementally in the classical framework). The new methods will potentially allow for significantly larger step sizes, similar to those reported in \cite{2018arXiv180106016M}, thereby substantially reducing the solution times.

We believe there are many opportunities of extending this approach to other problems and end our discussion by mentioning a few potential avenues. For extending the EN method to vector functions, there are other options as well, which need to be explored in the future. Similarly, simplifications of CN method to give quasi-Corrected Newton methods (or quasi-Halley's methods) would create a new area of research. In general, the extension of such root-finding methods to multiple-variables is a broad area, which also takes into account additional factors, such as sparsity and symmetry of the system, into account. The idea of multiplying with a function of $x$ instead of a constant, can be used in the development of preconditioners. Unlike linear preconditioners, which are well-understood, this will lead to nonlinear preconditioning, which is a relatively under-explored area \cite{cai2002nonlinearly}.

Application of a similar approach to nonlinear least squares fitting (or regression) will likely improve the fitting procedure, which is widely used in parameter estimation. Similarly, in optimization problems with multiple unknowns, a line search is commonly used. This approach could be used to improve the convergence of line search, which is akin to solving a single unknown equation. There is also the possibility of extending similar ideas to solving nonlinear dynamics, which is usually done by numerical time integration schemes. It is well-known that an accurate and stable solution of stiff ordinary differential equations, owing to their nonlinearity and coupling of different time scales, requires extremely small time steps. A similar approach of reducing their nonlinearity may be used to improve the time integration methods so that larger time steps can be used. 

It took several decades to develop methods based on the classical Newton's method. We anticipate that the new approach presented in this work, which resulted in some known (through classical higher-order approaches) and some new results, will spark the interest of the scientific community to fully explore the properties of these new methods, their applications, and the opportunities they present in science and engineering, over the next years.

\section*{Author Contributions}{A.A. and S.P. conceived the idea and designed the study. A.A. developed the theory and analysis. A.A. and S.P. analyzed the results and wrote the paper.}
\section*{Conflict of interest}{The authors declare no conflict of interest.}

\begin{acknowledgments}
{This work was supported by Welsh Government and Higher Education Funding Council for Wales through the S\^{e}r Cymru National Research Network in Advanced Engineering and Materials (Grant No. F28 to A.A.), and the Engineering and Physical Sciences Research Council of the UK (Grant No. EP/P018912/1 to A.A.).}
\end{acknowledgments}

\bibliography{refs}

\end{document}